\documentstyle[12pt, subeqnarray]{article}

\begin{document}
\newcommand{\dis}{\displaystyle}
\newcommand{\dml}{{\rm dim ~ ~ ker ~ ~}}
\newcommand{\expon}{{\rm e}}
\newcommand{\id}{ 1 \hspace{-2.85pt} {\rm I} \hspace{2.5mm}}
\newcommand{\beq}{\begin{equation}}
\newcommand{\eeq}{\end{equation}}
\newcommand{\su}{{\rm su}}
\newcommand{\U}{{\cal U}}
\newcommand{\Pa}{{\cal P}}
\newcommand{\bseq}{\begin{subeqnarray}}
\newcommand{\eseq}{\end{subeqnarray}}
\newcommand{\spa}{\mbox{\hspace{1cm}}}
\newcommand{\pr}{\prime}
{\thispagestyle{empty}
\rightline{} 
\rightline{} 
\rightline{} 

\vskip 0cm
\centerline{\large \bf Symmetric multiplets in}
\centerline{\large \bf Quantum Algebras}


\vskip 0.5cm
\centerline{
L. C. Kwek {\footnote{E-mail address: scip3057@nus.sg}}, 
C. H. Oh {\footnote{E-mail address:
phyohch@nus.sg } and K. Singh {\footnote{E-mail address:
crssingh@nus.sg }} } 
}
\centerline{{\it Department of Physics, Faculty of Science, } }
\centerline{{\it National University of Singapore,Lower Kent Ridge,} } 
\centerline{{\it Singapore 119260, Republic of Singapore. } }}
\vskip 0.1in


\vskip 3cm

\centerline{\bf Abstract} \vspace{10mm}
\noindent{
We consider a modified version of the coproduct for $\U(\su_q(2))$
and show
that in the limit when $q \rightarrow 1$, there exists  an
essentially non-cocommutative coproduct.
We study the implications of this non-cocommutativity
for a system of two  spin-$1/2$ particles.  Here it is shown that,
unlike the usual case, this non-trivial coproduct allows for symmetric
and anti-symmetric states to be present in the multiplet.  We surmise that
our analysis could be related to the
ferromagnetic and antiferromagnetic cases of the Heisenberg magnets.
 }

\newpage


In the standard Drinfeld-Jimbo\cite{drin} $ q$-deformation of
$\su(2)$ algebra which is defined by the relations
\beq
\lbrack J_+, J_- \rbrack = [2J_0] \mbox{\hspace{1cm}} [J_0 , J_\pm] = \pm
J_\pm \label{qcomm}
\eeq
\noindent where $\dis [x] = \frac{q^x - q^{-x}}{q - q^{-1}}$ and
 the coproduct is given by
\bseq
\Delta({J_{\pm}}) & =&  J_{\pm}\otimes q^{ J_0}
+ q^{- J_0}\otimes {J_{\pm}} \\
\Delta(J_0) & = &  J_0\otimes {\bf l}
+ {\bf l}\otimes J_0 , \label{qcocom}
\eseq
which is non-cocommutative for generic $q$, {\it i.e.} 
$\Delta (J_{\pm})\neq  \sigma \circ \Delta
(J_{\pm})$ and
$\sigma (a\otimes b)=b\otimes a$ is the flip automorphism. 
The coproduct dictates the tensor multiplications of two
representations. For instance in the case of two spin-$1/2$
representations the above coproduct leads to the following states:
\bseq
|0,0>& = &\frac{1}{\sqrt 2}(q^{\frac{1}{2}}|+>\otimes |->
-  q^{-\frac{1}{2}}|->\otimes | +>) \\
|1,1> & = & |+> \otimes |+> \\
|1,0> & =  & \frac {1}{\sqrt 2}(q^{-\frac{1}{2}}|+>\otimes | ->
+  q^{\frac{1}{2}}| ->\otimes |+>)\\
 | 1,-1> &=&  |->\otimes |->. \label{states}
\eseq

In ref[2], Zachos considered the $q \rightarrow -1$ limit of
this $\su_q(2)$ algebra and and showed that it has some interesting
consequences for the wavefunctions. Specifically, he showed
that for a system of two spin-1/2
particles, the singlet state which is ordinarily antisymmetric
transforms into a symmetric state while one of the triplet states becomes 
antisymmetric. This odd behaviour can be traced to the fact that
the coproduct under the $q \to -1$ limit remains noncocommutative
while the half-integer spin representations of $\su_q(2)$ reduce to
those of $\su(2)$.\footnote[4]{In the case of integer-spin representations, the
$q\to -1$ limit reduces to those of su(1,1).}
So on the surface, it appears that one has an unconventional
composition law for the usual su(2) algebra.
There is a caveat in
this argument however, the limit is singular in at least two respects.
Firstly, the coproduct obtained by taking the $q\to -1$ limit in (2a)
is not an algebra homomorphism  in the case of half integral spin
representations.\renewcommand{\thefootnote}{\fnsymbol{footnote}}{\protect
\footnote[2]{One can ascertain this by computing explicitly the
matrices of $\Delta J_\pm$ and $\Delta J_0$ in the spin $1/2 \otimes 1/2$
representation of $\su_q(2)$. By taking the limit $q \rightarrow -1$ of
these matrices one then finds that the commutator between $\Delta J_+$ and
$\Delta J_-$ yields $- 2 \Delta J_0$ instead of the requisite
$2 \Delta J_0$.}}
Secondly, as pointed out in ref[2] the Casmir becomes divergent
when one considers half integral spin representations:
\beq
\lbrack j \rbrack  \lbrack j + 1 \rbrack
\stackrel{q\rightarrow -1}{\longrightarrow}
\dis \frac{4}{\epsilon^2} + \frac{1}{2} + j(j + 1) 
~  + \{ -\frac{1}{32} + \frac{1}{24} (j(j + 1)(2j^2 +
2j - 1) \}\epsilon^2  + o(\epsilon^4)
\eeq
\noindent where $\epsilon = q - q^{-1}$.

In this letter, we look at the admissability of such states and ask
whether a non-trivial coproduct exists for the su(2) algebra. The latter
would  essentially lead to a universal enveloping algebra (UEA) of ordinary
$\su(2)$ 
Lie algebra with a non-cocommutative Hopf structure. We start by
modifying the $\su_q(2)$ algebra coproduct as:
\bseq
\Delta(J_+) & = & J_+ \otimes q^{J_0} \expon^{i n \pi J_0 } + \expon^{-
i n \pi J_0} q^{-J_0} \otimes J_+ \\
\Delta(J_-) & = & J_- \otimes q^{J_0} \expon^{- i n \pi J_0 } + \expon^{
i n \pi J_0} q^{-J_0} \otimes J_- \\
\Delta(J_0) & = & J_0 \otimes {\bf 1} + {\bf 1} \otimes J_0  \\
\epsilon(J_+) =  &  \epsilon(J_-) =  &  \epsilon(J_0) = 0 \\
S(J_\pm) = & -  \expon^{i n \pi} q^{\pm 1} J_{\pm},   &  S(J_0) = - J_0 \label{coprod}
\eseq
\noindent where $n$ is any integer.
In the limit when $q \rightarrow 1$,
 the coproduct, counit and antipode become
\bseq
\Delta(J_{\pm}) & = & J_{\pm}\otimes \expon^{\pm i\pi nJ_0}
+ \expon^{ \mp i\pi nJ_0}\otimes {J_{\pm}} \\
\Delta(J_0) & = &  J_0 \otimes {\bf l} + {\bf l}\otimes J_0  \\
\epsilon ({J_{\pm}}) & = & \epsilon (J_0) = 0 \\
S(J_{\pm}) & = & - \expon^{i\pi n}J_{\pm}\qquad S(J_0)=-J_0. \label{hopf1}
\eseq

It is important to observe that for {\bf odd $n$}, the coproduct
eq(\ref{hopf1} a) remains non-cocommutative.
To elucidate this non-triviality, it
is instructive to evaluate the coproduct and its opposite
$(\Delta'\equiv \sigma
\circ \Delta$) in a tensor product representation. To this end we
consider the tensor product of two representations labelled by $j_1$ and 
$j_2$ {\it i.e.} $(j_1\otimes j_2)$. Since the commutation relations 
are just the standard ones, we have for each representation
\bseq
J_0 |j,m> &= & m |j,m>\\
J_{\pm} |j,m> & = & \sqrt{(j\mp m)(j\pm m+1)} |j,m\pm 1> \label{stand}
\eseq
\noindent where $-j \leq m \leq j$.
The value $j$ characterizes the representation:
\beq
{\bf {J}}^2  |j,m> \equiv (J_+J_-+J_0^2- J_0) |j,m> = j(j+1) |j,m>
\eeq
\noindent where ${\bf J}^2$ is the quadratic casmir and $j$
takes the values $0,1/2,1,3/2, \cdots $. For an arbitrary vector
$|j_1,m_1>\otimes |j_2,m_2>$ belonging to the representation space of
$(j_1\otimes j_2)$ and using eq(\ref{hopf1}a), we have
\begin{eqnarray}
\Delta (J_{\pm}) |j_1,m_1>\otimes |j_2,m_2> & = &
f_{\pm}(j_1,m_1)\expon^{\pm i\pi nm_2} |j_1,m_1+1>\otimes |j_2,m_2>
\nonumber \\
& & + f_{\pm}(j_2,m_2)\expon^{\mp i\pi nm_1} |j_1,m_1>\otimes
|j_2,m_2+1> \label{rep1}
\end{eqnarray}
\noindent where $f_{\pm}(j_i,m_i)=
\sqrt{(j_i\mp m_i)(j_i\pm m_i+1)}$  and 
$-j_i\le m_i \le j_i$ ($i=1,2$). For the opposite coproduct one obtains
\begin{eqnarray}
\Delta'(J_{\pm})|j_1,m_1>\otimes |j_2,m_2> & = &
f_{\pm}(j_1,m_1)\expon^{\mp i\pi nm_2} |j_1,m_1+1>\otimes |j_2,m_2>
\nonumber\\
& & + f_{\pm}(j_2,m_2) \expon^{\pm i\pi nm_1} |j_1,m_1>\otimes
|j_2,m_2+1>. \label{rep2}
\end{eqnarray}

Now if the coproduct is cocommutative at the algebraic level it
is necessary that this property be reflected by any tensor product 
representation. In other words the rhs. of eq(\ref{rep1}) and
eq(\ref{rep2}) 
should be equal for {\it any} allowed 
values of $j_1$ and $j_2$ and their corresponding $m$'s. 
On the other hand, the {\it existence} of a
tensor product representation in which the two do not agree is
sufficient proof for non-cocommutativity of $\Delta$.  When one of the
two representations that appears in the tensor product is characterized
by a half-integer value ($ i.e. j_i = $ 1/2, 3/2, 5/2, ...), then for
odd integer $n$, we see that  
$\Delta (J_{\pm})\neq \Delta'(J_{\pm})$. So we
can surmise that $\Delta (J_{\pm})\ne \Delta'(J_{\pm})$ in general.

We stress that in the $q \rightarrow 1$ limit of eq(\ref{coprod}a,b),
the coproduct in (\ref{hopf1}a) differs from the usual one for
$\su_q(2)$ for $q = \expon^{in\pi}$.
By identifying $\dis \expon^{i \pi n} \equiv q^\prime$ in
(\ref{hopf1}a) one has
\beq
\Delta(J_\pm) = J_\pm \otimes q^{\prime\pm J_0} + q^{\prime\mp J_0}
\otimes J_\pm
\eeq
\noindent which clearly differs from eq(\ref{qcocom}a) by a sign in
$J_0$ for $\Delta(J_-)$.  Note further that this is a Hopf
$\ast$-algebra in the sense that the canonical conjugation\cite{singh}
\beq
(J_0)^+ = J_0 \spa (J_\pm)^\dagger = J_\mp
\eeq
\noindent is compatible with the Hopf structure.

For generic q, the tensor product of two representations closely
imitates the results in eq(\ref{rep1}), with the expressions
$f_\pm(j_i, m_i) $ replaced by its q-deformed form, $f_\pm^\prime$, where
\beq
f_{\pm}^\prime(j_i,m_i)=
\sqrt{[j_i\mp m_i][j_i\pm m_i+1]}.
\eeq

We next consider the implications of these results
for a system of two spin-1/2 particles. It is convenient to
introduce a new set of generators in which the
coproduct assumes an equivalent form. We define
\beq
J_+^\pr=\expon^{i\pi nJ_0}J_+,\qquad J_-^\pr=J_-\expon^{-i\pi
nJ_0}\qquad J_0^\pr=J_0
\eeq 
with the primed generators satisfying the same commutation relations
as the unprimed ones so that we can regard these
operators as the generators of original algebra.  The coproduct in
eq(\ref{coprod}a -c) now reads
\bseq
\Delta(J_+^\pr) & = & J_+^\pr \otimes q^{J_0^\pr} \expon^{2 i n \pi
J_0^\pr } + q^{-J_0^\pr} \otimes J_+^\pr \\
\Delta(J_-^\pr) & = & J_-^\pr \otimes q^{J_0^\pr} \expon^{- 2 i n \pi
J_0^\pr } + q^{-J_0^\pr} \otimes J_-^\pr \\
\Delta(J_0^\pr) & = & J_0^\pr \otimes {\bf 1} + {\bf 1} \otimes J_0^\pr. \label{newco}
\eseq
\noindent With $\dis j=\frac{1}{2}$ for spin-$1/2$ representation, we
have
\beq
J_{\pm}^\pr |\pm> = 0 \spa J_{\pm}^\pr |\mp> = |\pm> \spa J_0^\pr |\pm>
= \pm \frac{1}{2} |\pm>
\eeq
\noindent where $\dis |\pm>$ denotes the q-deformed states $\dis |\frac{1}{2},\pm
\frac{1}{2}>$ respectively.  Using the same techniques for evaluating the
normal Clebsch-Gordon coefficients based on the coproduct in
eq({\ref{newco}), the explicit expressions for the singlet and
triplet states
can be written as:
\bseq
|0,0>& = &\frac{1}{\sqrt 2}(q^{\frac{1}{2}}|+>\otimes |->
- (-1)^n q^{-\frac{1}{2}}|->\otimes | +>) \\
|1,1> & = & |+> \otimes |+> \\
|1,0> & =  & \frac {1}{\sqrt 2}(q^{-\frac{1}{2}}|+>\otimes | ->
+ (-1)^n q^{\frac{1}{2}}| ->\otimes |+>)\\
 | 1,-1> &=& (-1)^n |->\otimes |->. \label{qstates}
\eseq
\noindent The main difference from those resulting from the coproduct
(\ref{qcocom}a, b) is the appearance of the factor
$(-1)^n$. For even $n$, we get the usual q-deformed multiplets and in
the limit $q \rightarrow 1$, we retrieve the usual antisymmetric
singlet state and the usual symmetric $|1, 0>$ state in the triplet. For
odd $n$, in the limit $q \rightarrow 1$, we see that the
singlet state becomes symmetric while the $|1,0>$
state in the triplet becomes antisymmetric.
The states here are similar to those obtained through the $q \rightarrow
-1$ limit shown in ref[2].

In the spin-$1/2$ representation
eq (\ref{qstates}a - d), if we let 
$\Delta_0$ be the coproduct for even $n$ and $\Delta_1$ be the
coproduct for odd $n$, then the two coproducts are related by
$U \Delta_1(a) U^\dagger =
\Delta_0(a),$ where $a$ is $J_\pm^\prime$ or $J_0^\prime$ and $U$ is
given by
\beq
U = \left( \begin{tabular}{cccc}
1 & 0 & 0 & 0 \\
0 & 0 & -1 & 0 \\
0 & 1 & 0 & 0 \\
0 & 0 & 0 & -1
\end{tabular}
\right).
\eeq
\noindent However, we saw earlier that the representations for the two
coproducts are different and $\sigma \circ \Delta_1 \neq \Delta_1$ where
$\sigma$ is the usual flip homomorphism.
Now under this unitary transformation, the canonical exchange operator,
\beq
P  =  \left( \begin{tabular}{cccc}
1 & 0 & 0 & 0 \\
0 & 0 & 1 & 0 \\
0 & 1 & 0 & 0 \\
0 & 0 & 0 & 1
\end{tabular} \right),
\eeq
transform into
\beq
P ~ \rightarrow ~ P^\prime = ~  U P U^{\dagger}
~  =  ~ \left( \begin{tabular}{cccc}
1 & 0 & 0 & 0 \\
0 & 0 & -1 & 0 \\
0 & -1 & 0 & 0 \\
0 & 0 & 0 & 1
\end{tabular} \right),
\eeq
\noindent Such a redefinition is reminicent
of an analogous situation for
the Heisenberg magnet\cite{negele}
in which the exchange operator $P$ is defined by
\beq
P  =  \left( \begin{tabular}{cccc}
1 & 0 & 0 & 0 \\
0 & 0 & $\eta$ & 0 \\
0 & $\eta$ & 0 & 0 \\
0 & 0 & 0 & 1
\end{tabular} \right)
\eeq
\noindent where $\eta =1$ for ferromagnetic and $\eta = -1 $ for
antiferromagnetic cases.  Quantum mechanically, we also note that the
two cases (ferromagnetic and antiferromagnetic) are fundamentally
different.

To summarize briefly, we have shown that the $q \rightarrow 1$ limit of
the coproduct (\ref{hopf1}a,b) can lead to the
UEA of the su(2) algebra  endowed with a non-cocommutative Hopf
structure.  It is interesting to note that the
non-cocommutativity found here for su(2) does not have a
deformation parameter and its Hopf structure is fixed by
a particular choice of $n$. In this regard it differs from the
usual quantum algebras.
A consequence of this non-trivial coproduct is that it
does not commute with the flip operator $\sigma$, thus allowing
for symmetric and antisymmetric states to coexist within a muliplet.
In fact the existence of two coproducts reminds us
of a similar situation in
the Heisenberg magnets in which the exchange operator is defined
differently in the ferromagnetic and antiferromagnetic cases.
In passing we would like to remark that a similar Hopf structure
can be shown to exist for the UEA of su(3).

\end{document}